\documentclass[a4paper]{article}
\usepackage[latin1]{inputenc}
\usepackage[english]{babel}
\usepackage{amssymb}
\usepackage{amsmath}
\usepackage{amsthm}
\usepackage{mathrsfs}
\usepackage[all]{xy}
\usepackage[dvips]{graphicx}
\usepackage{url}
\newtheorem{teo}{Theorem}
\newtheorem{pro}[teo]{Proposition}
\newtheorem{lem}[teo]{Lemma}
\newtheorem{cor}[teo]{Corollary}
\newtheorem{rem}[teo]{Remark}
\newtheorem*{defi}{Definition}
\newtheorem*{con}{Conjecture}

\DeclareMathOperator{\Fix}{Fix}
\DeclareMathOperator{\id}{id}
\DeclareMathOperator{\Def}{Def}
\DeclareMathOperator{\Spec}{Spec}
\DeclareMathOperator{\coker}{Coker}
\newcommand{\of}{\mathcal{O}}

\newcommand{\pn}{\mathbb{P}}

\begin{document}
\title{\textbf{Symplectic involutions of holomorphic symplectic fourfolds}}
\author{Chiara Camere}
\date{ }
\maketitle
\section{Introduction}

In this paper we are going to study involutions of irreducible symplectic fourfolds and their fixed points. In particular we are going to concentrate on symplectic involutions, i.e. those which preserve the symplectic form.

The study of symplectic involutions and more generally of automorphisms of finite order on K3 surfaces has been started by Nikulin in \cite{Nik}. Since irreducible holomorphic symplectic manifolds are the natural generalization of K3 surfaces in higher dimension, Beauville started to study the same problems for such manifolds in \cite{Be3}. Many authors have studied the problem from different view-points, here we want to mention only the papers by Boissi\`{e}re \cite{Boi} and Boissi\`{e}re-Sarti \cite{Boi-Sar} on natural involutions and the paper by Beauville \cite{Be5} in which he deals with the case of antisymplectic involutions.

We remark in Section \ref{secfix} that the irreducible components of the fixed locus of a symplectic involution on an irreducible holomorphic symplectic manifold $ X $ are smooth symplectic subvarieties of $ X $; hence, if $ \dim X=4 $ they are either isolated points, smooth abelian surfaces or smooth K3 surfaces.

In Section \ref{sympl} we state and prove the main result of this paper: when $ X $ is an irreducible holomorphic fourfold with $ b_2(X)=23 $ there are only 3 possibilities for the number of isolated fixed points of a symplectic involution, 12, 28 or 36; moreover, in the first and third case at least one irreducible component of the fixed locus should be a smooth abelian surface. The main tool for our proof is the holomorphic Lefschetz theorem by Atyah-Singer (see \cite{AS}) that we recall in Section \ref{HLRR}.

In fact we conjecture that the fixed locus of a symplectic involution on such a fourfold cannot contain an abelian surface, and hence the only possibility would be that it contains 28 isolated points and a smooth K3 surface. The rest of the paper is devoted to provide evidences for this conjecture looking at different examples: the Hilbert scheme of a K3 surface in Section \ref{schema}, the Fano variety of a cubic in $ \pn^5 $ in Section \ref{Fano} and the double cover of an EPW sextic in Section \ref{lastOG} .

\section{Irreducible holomorphic symplectic manifolds}\label{IHS}

Let us recall first of all the definition of irreducible holomorphic symplectic manifold; for all the details on this subject the reader may refer to \cite{Be2} and to Part III of \cite{Huy1}.

\begin{defi}
A compact K\"{a}hler manifold $ X $ is irreducible holomorphic symplectic if it is simply connected and admits a symplectic 2-form $ \omega \in H^{2,0}(X) $ everywhere non degenerate and unique up to multiplication by a nonzero scalar.
\end{defi} 

It follows immediately from the definition that we have $ H^0(\Omega _X)\cong H^1(\of _X)=0 $, since $ X $ is simply connected; moreover the existence of a symplectic 2-form implies that the complex dimension of $ X $ is always even and that $ K_X $ is trivial. From the definition it follows that the Hodge structure of the second cohomology ring $ H^2(X,\mathbb{C})$ is $ H^2(X,\mathbb{C})\cong \mathbb{C}\omega\oplus H^{1,1}(X)\oplus\mathbb{C}\bar{\omega} $ and we have an isomorphism between $ TX $ and $ \Omega _X^1 $.

Not many examples of irreducible symplectic manifolds are known. Here we briefly describe those in dimension $ 4 $ that we need in the next sections.

\textbf{The Hilbert scheme of a K3 surface.} Let $ S $ be a smooth K3 surface and let $ X=S^{\left[2\right]} $ be the Hilbert scheme of $ S $ of $ 0- $schemes of length 2; then $ X $ can be constructed in the following way:

\[ \xymatrix{
\text{Bl}_{\Delta}(S\times S)\ar@{->}[r]\ar@{->}[d]&X\ar@{->}[d]\\
S\times S\ar@{->}[r]& S^{(2)} } \]
as the blow-up along the diagonal $ \Delta $ of the symmetric product of $ S $. For further details see \cite{Be2}. All the other families we are going to consider turn out to be deformation equivalent to this one and hence will have the same cohomology.

In particular let us recall that $ b_2(X)=23 $ and that the Hodge diamond is the following
\begin{center}
\begin{tabular}{ccccccccc}
&&&&1&&&&\\
&&&0&\ &0&&&\\
&&1&\ &21&\ &1&&\\
&0&\ &0&\ &0&\ &0&\\
1&\ &21&\ &232&\ &21&\ &1\\
\end{tabular}
\end{center}

\textbf{The Fano variety of a smooth cubic in $\pn^5$.} Let $ X $ be a smooth cubic hypersurface in $ \pn^5 $ and $ F $ the Fano Variety of $ X $, i.e. the variety of lines contained in $ X $. In their paper of 1985 \cite{B-D} Beauville and Donagi show that this is an irreducible holomorphic symplectic fourfold deformation equivalent to the former family. Moreover, they show that there is an isomorphism of Hodge structures $\alpha :H^4(X,\mathbb{Z})\cong H^2(F,\mathbb{Z})$.

\textbf{The double cover of an EPW sextic.} This example has been introduced and intensively  studied by O'Grady in \cite{OG} and many other papers.  Starting from  a 6-dimensional vector space $ V $  and from a general enough Lagrangian subspace $ A \subset \wedge^3V$, the subvariety $Y_A:=\left\lbrace v\in \pn (V)/(v\wedge\wedge^2 V)\cap A\neq 0\right\rbrace  $ turns out to be a hypersurface  of degree 6 of the type described by Eisenbud-Popescu-Walter; $ Y_A $ is not smooth, but it has a smooth double cover $ X_A $ that is an irreducible holomorphic symplectic fourfold when $ A $ is general enough. We will discuss this more in detail in Section \ref{lastOG}.

\section{Holomorphic Lefschetz Theorem}\label{HLRR}

Let us briefly recall the Holomorphic Lefschetz theorem by Atiyah-Singer (see \cite{AS}), following the paper by Donovan \cite{Don}, where the reader can find all the details and the proofs that we are skipping. In order to keep the notation as simple as possible we limit our presentation to the case of involutions.


Let $ Z\subset \Fix(i) $ be an irreducible component of the fixed point set of an involution $ i $ on a smooth projective variety $ X $. Let $ N_Z^* $ be the dual of the normal bundle of $ Z $; since $ TZ $ is fixed by $ di $ and on the other hand $ i $ is non degenerate, from the exact sequence
\[ \xymatrix{
0\ar@{->}[r]& TZ\ar@{->}[r]&TX _{|Z}\ar@{->}[r]&N_Z\ar@{->}[r]&0 } \]
it follows that $ N_Z $ is the eigensheaf corresponding to $ -1 $.

Let us consider a coherent sheaf $ F $ on $ X $ and let $ \eta:i^*F\longrightarrow F $ be a morphism such that the composite morphism 
$\eta\circ i^*\eta$ is the identity; then there is an induced action $ i^* $
on the vector space of global sections $ \Gamma(F) $. Hence the involution $ i $ and $ \eta $ induce an action on the cohomology ring $ H^*(X,F) $ of the sheaf $ F $ that we will always denote by $ i^* $ for the sake of simplicity.

Moreover $ \eta $ induces an involution $ \eta _{|Z}: F_{|Z}\longrightarrow F_{|Z} $ with eigenvalues $ \pm 1 $. Throughout all that follows we will denote $F^+_{|Z} $ and $ F^-_{|Z} $ the eigensheaves respectively fixed by $ \eta _{|Z} $ and associated to $ -1 $.


\begin{teo}\label{lrrf}
\textbf{Holomorphic Lefschetz-Riemann-Roch formula} Let $ X $ be a smooth projective variety of dimension $ d $, $ i $ an involution of $ X $ and $ F $ a coherent sheaf on $ X $; let $ \eta:i^*F\longrightarrow F $ be a morphism such that $\eta\circ i^*\eta=\id _F$ and let $ i^* $ be the induced action on $ H^*(X,F) $. Let $ Z\subset \Fix(i) $ be an irreducible component of the fixed point set. Then the following formula holds:
\[ 
\sum _{j=0}^d(-1)^j \mathrm{Tr}(i^*_{|H^j(F)})=\sum _{^{Z\subset \Fix(i)} _{irreducible}}\int _Z\dfrac{\mathrm{Td}(Z).\left[\mathrm{ch}(F_{|Z}^+)-\mathrm{ch}(F_{|Z}^-)\right]}{\sum_{p\geq 0} \mathrm{ch}(\wedge ^p N_{Z}^*)}
  \]
\end{teo}

\section{Fixed loci}\label{secfix}

The property of preserving the symplectic form induces limitations on the irreducible components of the locus of fixed points. Let us remark some important properties.

\begin{lem}\label{lsmo}
Let $ X $ be a projective smooth variety and $ f:X\longrightarrow X $ a periodic endomorphism; then each component of the fixed point set $ \Fix(f) $ is smooth.
\end{lem}
\textbf{Proof.} See \cite{Don}, Lemma 4.1. \qed\\

In the case we are interested in we have more than smoothness.

\begin{pro}\label{psymsub}
Let $ X $ be an irreducible holomorphic symplectic manifold and $ i $ a symplectic involution on $ X $.
Then the irreducible components of $ \Fix(i) $ are symplectic subvarieties of $ X $.
\end{pro}

\textbf{Proof.} Let $ Z\subset \Fix(i) $ be an irreducible component of the fixed point set of dimension $ d>0 $. We need to prove that the restriction to $ Z $ of the symplectic form $ \omega $ gives a symplectic form on $ Z $. We know that $ TX _{|Z}\cong TZ\oplus N_Z $ and that $ TZ $ and $ N_Z $ are respectively the eigensheaves associated to $ \pm 1 $. Given $ z\in Z $, since $ i  $ is symplectic, $ T_zZ $ and $ N_{Z,z} $ are orthogonal and hence both symplectic. \qed


\begin{rem}
In particular if $ X $ has dimension 4, the irreducible components can be either isolated fixed points or K3 and abelian surfaces.
\end{rem}

\section{Symplectic involutions}\label{sympl}

Now we are ready to show the main result. We will show that there are few different possibilities for the nature of the fixed locus of a symplectic involution on an irreducible holomorphic symplectic fourfold such that $ b_2=23 $. By \cite{G}, when $X$ is an irreducible holomorphic symplectic fourfold such that $ b_2=23 $, the canonical map $S^2H^2(X,\mathbb{C})\rightarrow H^4(X,\mathbb{C})$ is an isomorphism and, as we already said in Section \ref{IHS}, this is the case for the family of Hilbert schemes $S^{\left[2\right]} $ of a K3 surface $S$ and for their deformations.

\begin{teo}\label{syminv}
Let $ X $ be an irreducible holomorphic symplectic fourfold such that $ b_2(X)=23 $ and let $ i $ be a symplectic involution of $ X $. Let $ \tau $ be the trace of $ i^* $ on $ H^{1,1}(X) $, $ N $ and $ K $ respectively the numbers of isolated fixed points and of K3 surfaces of fixed points. Then only the following cases can occur:
\begin{enumerate}
\item $\tau=-3$, $N=12$ and $K=0$;
\item $\tau=3$, $N=36$ and $K=0$;
\item $\tau=5$, $N=28$ and $K=1$.
\end{enumerate}
Moreover in the first two cases $ i $ fixes at least one abelian surface.
\end{teo}

 In fact, we conjecture that only the last case can occur.

\begin{con}
 Let $ X $ and $ i $ be as in Theorem \ref{syminv}; the fixed locus of $ i $ cannot contain an abelian surface.
 \end{con}

In the next sections we will provide evidence for this conjecture verifying it in some of the known examples of irreducible symplectic fourfolds such that $ b_2(X)=23 $.\\

\textbf{Proof of Theorem \ref{syminv}.} Let us apply the holomorphic Lefschetz Riemann-Roch formula discussed in Section \ref{HLRR} to the cohomology of the sheaves $\of _X$, $\Omega ^1 _X$ and $\Omega ^2 _X$.

\textbf{The sheaf} $\mathbf{\of _X}.$ We know that $h^{0,0}=h^{2,0}=h^{4,0}=1$ and $h^{1,0}=h^{3,0}=0$; on the other hand, we know that $H^{2,0}(X)=\langle\omega _X\rangle$ and $H^{4,0}(X)=\langle\omega _X^2\rangle$, hence they are fixed by the involution and the Lefschetz number is $L(i)=\Sigma (-1)^i\mathrm{Tr}(i^*_{|H^{i,0}(X)})=3$.

For each fixed surface $Y$ of $\Fix\left(i\right) $ we have to calculate $$\int_Y \text{Td}(Y).( 1+ \text{ch}(N^*_Y)+\text{ch}(\det N^*_Y)) ^{-1}$$ since the only eigenvalue of $\text{d}i$ is -1 and the rank of $N^*_Y$ is 2. From the short exact sequence
\[ 
\xymatrix{
0\ar@{->}[r]&TY\ar@{->}[r]&TX_{|Y} \ar@{->}[r]&N_Y\ar@{->}[r]&0}\]
we get $c_1(N^*_Y)=-c_1(N_Y)=-c_1(TX_{|Y})+c_1(TY)=0$ and $$c_2(N^*_Y)=c_2(N_Y)=c_2(TX_{|Y})-c_2(TY)-c_1(N_Y)c_1(TY)=c_2(X).\left[ Y\right]-c_2(Y). $$

The Lefschetz formula for $\of _X$ becomes then
\begin{equation}\label{oxprima}
 3=\sum _{i(p)=p}\frac{1}{\det (1-\text{d}i_{|T_p})}+\sum _{i(S_j)=S_j}\int _{S_j}\frac{1+\frac{1}{12}c_2(S_j)}{4+c_2(S_j)-c_2(X).\left[S_j\right]}
\end{equation}

The restriction of $\text{d}i$ to $T_p$ is $-\id_{\mathbb{C}^4}$, hence $\det (1-\text{d}i_{|T_p})=2^4=16$.

An easy computation gives $$\frac{1+\frac{1}{12}c_2(S_j)}{4+c_2(S_j)-c_2(X).\left[S_j\right]}=\frac{1}{4}-\frac{1}{24}c_2(S_j)+\frac{1}{16}c_2(X).\left[S_j\right]$$

Let us write $a_j:=\int _{S_j}c_2(X).\left[S_j\right]$; since $\int _{S_j}c_2(S_j)=24$ if $S_j$ is K3 and $\int _{S_j}c_2(S_j)=0$ if $S_j$ is abelian, from (\ref{oxprima}) it follows that

\begin{equation}
 3=\frac{N}{16}-K+\frac{1}{16}\sum_{i(S_j)=S_j}a_j
\end{equation}

\textbf{The sheaf} $\mathbf{\Omega ^1 _X}.$
We know that $h^{0,1}=h^{2,1}=h^{4,1}=0$ and $h^{1,1}=h^{3,1}=21$, since $H^{1,1}(X)\cong H^{3,1}(X)$; moreover this isomorphism is compatible with $ i^* $ since it is given by product with $ \omega_X $ and $ i $ is symplectic.
The Lefschetz number is $L(i,\Omega ^1 _X)=-2\tau$.

In this case we get
$$
L(i,\Omega ^1 _X)=\sum _{i(Y)=Y}\int_Y \frac{\text{Td}(Y).(\text{ch}(\Omega ^+ _{|Y})-\text{ch}(\Omega ^- _{|Y}))}{( 1+ \text{ch}(N^*_Y)+\text{ch}(\det N^*_Y))}
$$
where $\Omega ^+ _{|Y}$ and $\Omega ^- _{|Y}$ are respectively the subbundle of $\Omega ^1 _{X|Y}$ fixed by the action of the dual of $\text{d}i_{|Y}$ and the subbundle on which the dual of $\text{d}i_{|Y}$ has eigenvalue -1.

When $Y$ is an isolated fixed point $p\in X$ we get $$\frac{-4}{\det (1-\text{d}i_{|T_p})}=-\frac{1}{4}$$

When $Y$ is a fixed surface, we have to calculate $(\text{ch}(\Omega ^+ _{|Y})-\text{ch}(\Omega ^- _{|Y}))$. Since $\Omega ^+ _{|Y}=TY$ and $\Omega ^- _{|Y}=N_Y^*$, we have
$$
(\text{ch}(\Omega ^+ _{|Y})-\text{ch}(\Omega ^- _{|Y}))=2-c_2(Y)-2+c_2(X).\left[ Y\right]-c_2(Y)=c_2(X).\left[ Y\right]-2c_2(Y)
$$
From $$ \left( \frac{1}{4}-\frac{1}{24}c_2(Y)+\frac{1}{16}c_2(X).\left[Y\right]\right) .\left( c_2(X).\left[ Y\right]-2c_2(Y)\right)=\frac{1}{4}c_2(X).\left[ Y\right]-\frac{1}{2}c_2(Y)
$$
it follows then that the Lefschetz formula for $\Omega ^1 _X$ is
\begin{equation}
-2\tau=-\frac{N}{4}-12K+\frac{1}{4}\sum_{i(S_j)=S_j}a_j
\end{equation}

\textbf{The sheaf} $\mathbf{\Omega ^2 _X}.$
We know that $h^{0,2}=h^{4,2}=1$, $h^{2,2}=232$ and $h^{1,2}=h^{3,2}=0$: indeed, from $S^2H^2(X,\mathbb{C})\cong H^4(X,\mathbb{C})$ (see \cite{G}) it follows $$H^{2,2}(X)\cong H^{2,0}(X)\otimes H^{0,2}(X)\oplus S^2H^{1,1}(X)\cong \mathbb{C}\oplus S^2H^{1,1}(X)$$ Let us write $\sigma:=\mathrm{Tr}(i^*_{|H^{2,2}(X)})$; we need to deduce $\sigma$ from $\tau$.

If $i^*$ is of type $(a,b)$ on $H^{1,1}(X)$, we have $\tau=a-b$ and $h^{1,1}=a+b$; on the other hand $H^{1,1}(X)=H^{1,1}_+\oplus H^{1,1}_-$ implies that $$H^{2,2}(X)\cong \mathbb{C}\oplus S^2H^{1,1}_+\oplus S^2H^{1,1}_-\oplus H^{1,1}_+\otimes H^{1,1}_-$$
Hence $\sigma=1+\frac{a(a+1)}{2}+\frac{b(b+1)}{2}-ab=1+\frac{h^{1,1}}{2}+\frac{\tau ^2}{2}$.

The Lefschetz number is then $L(i,\Omega ^2 _X)=2+\sigma=3+\frac{21}{2}+\frac{\tau ^2}{2}=\frac{27+\tau ^2}{2}$.

We need to know which are the subbundles of $\Omega ^2 _{X|Y}$ associated to the eigenvalues 1 and -1. We have $(\Omega ^2 _{X|Y})^+\cong \wedge ^2\Omega ^+ _{|Y}\oplus\wedge ^2\Omega ^- _{|Y}\cong\det TY\oplus \det N^*_Y$ and $(\Omega ^2 _{X|Y})^-\cong TY\otimes N^*_Y $.

When $Y$ is an isolated fixed point $p\in X$ we get $$\frac{6}{\det (1-\text{d}i_{|T_p})}=\frac{3}{8}$$

When $Y$ is a fixed surface, we get instead
$$
(\text{ch}(\Omega ^2 _{X|Y})^+)-\text{ch}(\Omega ^2 _{X|Y})^-))=\text{ch}(\of _Y^2)-\text{ch}(TY)\text{ch}(N^*_Y)=2-4+2c_2(X).\left[ Y\right]=$$ $$=2c_2(X).\left[ Y\right]-2
$$
Since $$ \left( \frac{1}{4}-\frac{1}{24}c_2(Y)+\frac{1}{16}c_2(X).\left[Y\right]\right) .\left( 2c_2(X).\left[ Y\right]-2\right)=-\frac{1}{2}+\frac{3}{8}c_2(X).\left[ Y\right]+\frac{1}{12}c_2(Y)
$$
the Lefschetz formula for $\Omega ^2 _X$ becomes
\begin{equation}
\frac{27+\tau ^2}{2}=\frac{3N}{8}+2K+\frac{3}{8}\sum_{i(S_j)=S_j}a_j
\end{equation}

\textbf{The system.}
We have thus obtained the following linear system
\begin{equation}\label{sys3}
\left\lbrace \begin{array}{l}
 3=\frac{N}{16}-K+\frac{1}{16}\sum_{i(S_j)=S_j}a_j\\
 \\
 -2\tau=-\frac{N}{4}-12K+\frac{1}{4}\sum_{i(S_j)=S_j}a_j\\
 \\
\frac{27+\tau ^2}{2}=\frac{3N}{8}+2K+\frac{3}{8}\sum_{i(S_j)=S_j}a_j
\end{array}
\right.
\end{equation}

from which, by eliminating $ \sum_{i(S_j)=S_j}a_j $, we deduce 

\begin{equation}\label{syssep}
\left\lbrace \begin{array}{l}
 -\tau ^2+4\tau +33=N\\
 \\
\tau ^2-9=16K
\end{array}
\right.
\end{equation}

On the other hand it must be $ K\geq 0 $ and $ N\geq 0 $, hence $ \tau $ must satisfy the following:
$$
\left\lbrace \begin{array}{l}
 \tau\leq -3\text{ or }\tau\geq 3 \\
 \\
2-\sqrt{37}\leq\tau\leq 2+\sqrt{37}
\end{array}
\right.
$$
Moreover the second equation of (\ref{syssep}) implies that $ \tau $ is odd and that it cannot be $ 7 $, since otherwise $ K $ would not be an integer. 

 If we replace the solutions we found in (\ref{sys3}), we get
\begin{enumerate}
\item $ \sum_{i(S_j)=S_j}a_j =36$ when $ \tau=-3 $;
\item $ \sum_{i(S_j)=S_j}a_j =12$ when $ \tau=3 $;
\item $ \sum_{i(S_j)=S_j}a_j =36$ when $ \tau=5 $.
\end{enumerate}
Hence if there is a symplectic involution satisfying the first or the second line of the table it must have a fixed abelian surface. This ends the proof.\qed

\begin{cor}
Let $ X $ and $ i $ be as in Theorem \ref{syminv}; then:
\begin{enumerate}
\item $ i $ has always at least 12 isolated fixed points and 1 fixed surface;
\item $ i $ fixes at most 1 K3 surface and in this case it has 28 isolated fixed points.
\end{enumerate}                 
\end{cor}

\section{The Hilbert scheme of a K3 surface}\label{schema}

As a first evidence to our conjecture, we will show that the natural symplectic involution on the Hilbert scheme of a K3 surface fixes exactly 28 isolated points and 1 K3 surface.

Let $ S $ be a smooth K3 surface and let $ X $ be the Hilbert scheme of $ S $ of $ 0- $schemes of length 2 (see Section \ref{IHS} for some details on the construction); given an involution $ \sigma $ of $ S $, there is an involution $ i=\sigma^{\left[ 2\right] } $ induced by it: such an involution is said to be natural. For further details on natural involutions the reader is referred to \cite{Boi-Sar} and \cite{Boi}.

Here we want to remark only that if $ \sigma $ is symplectic then also $ i $ will preserve the symplectic form on $ X $. Moreover, Nikulin showed in \cite{Nik} that a symplectic involution on a smooth K3 surface fixes 8 isolated points. Hence, on $ X $ the isolated fixed points will be all the couples $ \lbrace p,q\rbrace $ where $ p,q\in \Fix(\sigma) $ are distinct; this gives $\binom{8}{2}=28$ isolated fixed points. The fixed K3 surface is the closure in $ X $ of the surface made of the points $ \lbrace p,\sigma(p)\rbrace $ with $ p\in S\setminus \Fix(\sigma) $.

Let us study the deformations of the couple $ (X,i) $. We will show that there are nontrivial deformations, i.e. deformations that cannot be obtained from a deformation of $(S,\sigma)$.

 The infinitesimal deformations of $ X $ are unobstructed and there is a canonical isomorphism $ j :H^2(S,\mathbb{C})\oplus \mathbb{C}e \longrightarrow H^2(X,\mathbb{C})$ (see \cite{Be2}), where $ e $ is the class of the exceptional divisor.  


\begin{pro}\label{defhilb}
Let $ S $ be a smooth K3 surface and $ \sigma $ a symplectic involution on $ S $; let $ X=S^{\left[2\right]} $ be the Hilbert scheme of $ S $ and $ i=\sigma ^{\left[2\right]} $ the natural symplectic involution on $ X $. Then the infinitesimal deformations of the couple $ (X,i) $ are parametrized by $H^{1,1}(X)^{i}=j(H^{1,1}(S)^{\sigma})\oplus \mathbb{C}e$.
\end{pro}

\textbf{Proof.} We have the following diagram
\[ 
\xymatrix{
\Def(S,\sigma )\ar@{->}[r]\ar@{->}[d]&\Def(S)\ar@{->}[d]\\
\Def(X,i)\ar@{->}[r]&\Def(X)}\]

Looking at the tangent spaces at 0 we thus get
\[ 
\xymatrix{
H^{1,1}(S)^{\sigma}\ar@{->}[r]\ar@{->}[d]&H^{1,1}(S)\ar@{->}[d]\\
H^{1,1}(X)^{i}\ar@{->}[r]&H^{1,1}(X)=j(H^{1,1}(S))\oplus \mathbb{C}e}\]
Since $\tau=5$ and $h^{1,1}(S)=20$, we have $\dim H^{1,1}(X)^{i}=13$; on the other hand
$\dim H^{1,1}(S)^{\sigma }=12$ by Theorem \ref{HLRR} and all natural automorphisms leave globally invariant the exceptional divisor (see \cite{Boi-Sar}), hence $i^*e=e$ and $e\in H^{1,1}(X)^{i}$. As a consequence we see that $H^{1,1}(X)^{i}=j(H^{1,1}(S)^{\sigma})\oplus \mathbb{C}e$.\qed



\section{The Fano variety of a smooth cubic}\label{Fano}

Let $X$ be a smooth cubic in $\pn ^5$ and let $F$ be the variety of lines of $X$; it is an irreducible holomorphic symplectic fourfold (see \cite{B-D}). We want to investigate which involutions $\sigma$ of $\pn ^5$ induce involutions of $X$ and hence of $F$ and of which kind these ones are.

We have the following situation:

\begin{center}
\begin{tabular}{|c|c|c|}
\hline
$\sigma$ & $ X $ & $p\in X$ s.t.. $\sigma(p)=p$\\
\hline
$\left[X_0,..,X_5\right] $&$X_0^2L+G$&$\left[1,0,..,0\right]$,\\
$\downarrow$& with $L\in \mathbb{C}\left[X_1,..,X_5\right]_1$&$\left[0,y_1,..,y_5\right]\in V(G)$\\
$\left[-X_0,X_1,..,X_5\right]$& $G\in \mathbb{C}\left[X_1,..,X_5\right]_3$&\\
\hline
$\left[X_0,..,X_5\right]$ &$X_0^2L_0+X_1^2L_1+X_0X_1L_2+G$&$\left[x_1,x_2,0,..,0\right]$ $\forall\left[x_1,x_2 \right]\in \pn ^1 $,\\
$\downarrow$& with $L_i\in \mathbb{C}\left[X_2,..,X_5\right]_1$ &$\left[0,0,y_1,..,y_4\right]\in V(G)$\\
$\left[-X_0,-X_1,X_2,..,X_5\right]$& $G\in \mathbb{C}\left[X_2,..,X_5\right]_3$&\\
\hline
$\left[X_0,..,X_5\right]$ &$X_0^2L_0+\dots+X_2^2L_5+G$&$\left[x_1,x_2,x_3,0,0,0\right]$\\
$\downarrow$& with $L_i\in \mathbb{C}\left[X_3,..,X_5\right]_1$ & $\forall\left[x_1,x_2,x_3 \right]\in \pn ^2 $,\\
$ \left[-X_0,-X_1,-X_2,X_3,X_4,X_5\right]$& $G\in \mathbb{C}\left[X_3,..,X_5\right]_3$&$\left[0,0,0,y_1,..,y_3\right]\in V(G)$\\
\hline
\end{tabular}
\end{center}
\vspace{0.5 cm}
In \cite{B-D} the authors show that $\alpha :H^4(X,\mathbb{Z})\cong H^2(F,\mathbb{Z})$ is an isomorphism of Hodge structures; via this isomorphism we have $\alpha(H^{2,0}(F))=H^{3,1}(X)$. By Griffiths' theorem on the cohomology of hypersurfaces in $\pn ^n$ (see \cite{Voi} \S 18 Th\'eor\`eme 18.1), 
$H^{3,1}(X)$ is generated by a meromorphic 5-form of $\pn ^5$ with poles of order 2 along $X$, i.e. $$\Omega = \sum (-1)^i X_i\frac{dX_0\wedge..\wedge \hat{dX_i}\wedge..\wedge dX_5}{P^2} $$ where $P$ is a polynomial defining $X$. Hence $\sigma$ induces on $F$ a symplectic involution $ i $ if and only if $\sigma^*\Omega=\Omega$ and this is true only in the second case.

Let us study in more detail the locus of fixed points on $F$ in 
the  symplectic case: we claim that the fixed locus is given by 28 isolated points and one K3 surface.

Indeed there are the lines of fixed points, i.e. $ X_2=...=X_5=0 $ and the 27 lines on the cubic $ G=X_0=X_1=0 $, which give the 28 points. All other lines fixed by the involution pass through 2 fixed points, hence they can be parametrized as\[ \lambda \left[ a_1,a_2,0,..,0\right] +\mu\left[0,0,b_1,..,b_4 \right]  \] 
Replacing in the equation of $ X $ we get
\[ a_1^2L_0(b_1,..,b_4)+a_2^2 L_1(b_1,..,b_4)+a_1a_2L_2(b_1,..,b_4)=0\]
equation which defines a divisor of bidegree $ (2,1) $ in $ \pn^1\times V(G)\subset \pn^1\times \pn^3 $, which is a K3 surface.


When $i$ is symplectic, let us study $\Def(F,i)$ and compare it with $\Def(X,\sigma)$. 


\begin{pro}
Let $X$ be a smooth cubic in $\pn ^5$ and let $F$ be the variety of lines of $X$; let $\sigma$ be the involution of $\pn ^5$ such that $\sigma^*\Omega=\Omega$ and $i$ the symplectic involution induced by $\sigma$ on $F$. Then the infinitesimal deformations of the couple $(F,i)$ are parametrized by $ H^{1,1}(F)^i\cong H^1(X,TX)^{\sigma}$.
\end{pro}
\textbf{Proof.} We have (see \cite{Voi} Corollary 18.12 and Lemma 18.15) $H^1(X,TX)\cong R_P^3\cong H^{2,2}(X)_0$, where $R_P^3$ is the degree 3 component of the Jacobian ring of $P$, but on the other hand (see \cite{B-D}) we know that $H^{1,1}(F)_0\cong H^{2,2}(X)_0$. Hence also the invariant parts will be isomorphic, i.e. all deformations of $(F,i)$ are obtained by deforming $(X,\sigma)$ and taking the Fano variety of the deformation with the induced involution. \qed


\section{Double covers of EPW sextics}\label{lastOG}

As a last example let us see what happens in the case of the double cover of an EPW sextic. Let $ V $ be a 6-dimensional vector space, $ \pn(V)\cong\pn ^5 $; on $ \wedge^3V $ the wedge product $ \wedge:\wedge^3V\times\wedge^3V \longrightarrow \wedge^6V  $ induces a sympletic form $ \omega $ by choosing an isomorphism $ \wedge^6V\cong\mathbb{C} $. Let us take a Lagrangian subspace $ A \subset \wedge^3V$ and let us define $Y_A:=\left\lbrace v\in \pn (V)/(v\wedge\wedge^2 V)\cap A\neq 0\right\rbrace  $; for general $ A $, $ Y_A $ is a hypersurface  of degree 6 of the type described by Eisenbud-Popescu-Walter. Such a hypersurface is not smooth, but for a general $ A $ it has a smooth double cover $ X_A $ which is an irreducible holomorphic symplectic fourfold (see \cite{OG}). We will show that the involution of $V$ such that $\dim V^+=4$ induces, when $A$ is general enough, a symplectic involution on $X_A$ which fixes exactly 28 isolated points and one K3 surface.

Let $ F $ be the vector bundle given on fibers by $ F_v=v\wedge\wedge^2 V $ for all $ v\in \pn^5 $; there is an isomorphism $ F\cong\Omega _{\pn^5}^3(3) $. We look at the morphism $ \lambda _A :F\longrightarrow \frac{\wedge ^3 V}{A}\otimes \of _{\pn^5}$ given on fibers by $$ v\wedge \alpha\in F_v\mapsto \left[ v\wedge \alpha\right] \in \frac{\wedge ^3 V}{A}  $$ for $ v\in\pn^5 $. Since the two sheaves considered have both rank 10, $ \coker\lambda _A $ is a torsion sheaf on $ \pn^5 $. 

Let us define $\xi_A=\zeta_A\otimes \of_{Y_A}(-3)  $ where $ \zeta_A $ is a coherent sheaf on $ Y_A $ such that $ j_*\zeta_A=\coker\lambda _A $; there is an isomorphism $ \alpha_A:\xi _A\longrightarrow \xi _A^* $ that gives $ \of _{Y_A}\oplus\xi_A $ the structure of a commutative $ \of _{Y_A} -$algebra. We define $ X_A= \Spec (\of_{Y_A}\oplus\xi_A) $; the structure map $ f:X_A\longrightarrow Y_A $ is finite of degree 2. The fourfold
$ X_A $ is smooth whenever $$ A\in \mathbb{LG}^0(\wedge ^3 V)=\lbrace A\in  \mathbb{LG}(\wedge ^3 V)/\pn(A)\cap G(3,6)=\emptyset\ \mathrm{ and }\ \dim (v\wedge\wedge ^2 V)\cap A\leq 2\ \mathrm{ for\  all }\ v\in V\rbrace$$ The double cover is ramified over $ W_A=\lbrace v\in Y_A/\dim (v\wedge\wedge ^2 V)\cap A= 2\rbrace $. 
\begin{lem}
Let $ i:V\longrightarrow V $ be an involution; if $ i(A)=A $, $ i $ induces an involution $ \hat{\imath} $ on $ X_A $.
\end{lem}
\textbf{Proof.}
We clearly have $ i^*\Omega _{\pn^5}^3(3)\cong \Omega _{\pn^5}^3(3) $.
$$
\xymatrix{
F_{i(v)}\ar@{->}[r]\ar@{->}[d]^i& \frac{\wedge ^3 V}{i(A)} \ar@{->}[d]^i\\
F_v \ar@{->}[r]^{\lambda _{A,v}}& \frac{\wedge ^3 V}{A}}$$
and hence also
\begin{equation}\label{if}
\xymatrix{
i^*F\ar@{->}[r]\ar@{->}[d]& \frac{\wedge ^3 V}{i(A)}\otimes \of _{\pn^5} \ar@{->}[d]\\
F \ar@{->}[r]^{\lambda _A}& \frac{\wedge ^3 V}{A}\otimes \of _{\pn^5}}\end{equation}

If $i(A)=A$ then $ i(Y_A)=Y_A $. Indeed, $Y_A=\left\lbrace v\in \pn (V)/(v\wedge\wedge^2 V)\cap A\neq 0\right\rbrace  $ is invariant for $i$ as soon as $A$ is. 

In order to prove that $ i $ induces an involution on $ X_A $ we also need to show that $ i^*\zeta_A\cong\zeta_A $ and that the morphism $ \alpha_A:\zeta_A\longrightarrow\zeta_A^* $ commutes with $ i $. It follows from diagram (\ref{if}) that $ i^*\coker \lambda _A\cong \coker \lambda _A $ and this implies $ i^*\zeta_A\cong\zeta_A $.

Moreover we know that $ \alpha _A=\beta _A\otimes \id _{\of_{Y_A}(-3)} $ (see \cite{OG} Proof of Proposition 4.4) where $ \beta _A $ satisfies the following diagram
$$
\xymatrix{
0\ar@{->}[r]&F\ar@{->}[r]^{\lambda _A}\ar@{->}[d]& A^*\otimes \of _{\pn^5} \ar@{->}[r]\ar@{->}[d]& j_*\zeta_A\ar@{->}[r]\ar@{->}[d]^{\beta _A}&0\\
0\ar@{->}[r]&A\otimes \of _{\pn^5}\ar@{->}[r]^{\lambda _A^*}&F^*\ar@{->}[r]& \text{Ext}^1(j_{\ast}\zeta _A,\of_{\pn^5})\ar@{->}[r]&0}
$$
Everything in this diagram is invariant for $ i $, hence $ \beta_A $ and consequently $ \alpha_A $ commutes with $ i $.
This shows that $ X_A $ is fixed by an involution $ \hat{\imath} $ such that
\begin{equation}\label{comminv}
 \xymatrix{
X_A\ar@{->}[r]^{\hat{\imath}}\ar@{->}[d]^f& X_A\ar@{->}[d]^f\\
Y_A\ar@{->}[r]^{i}& Y_A}
 \end{equation}
\qed
\begin{rem}\label{indinv}
Given $ i $ involution on $ Y_A $, there are two involutions $ i_1 $ and $ i_2 $ on $X_A$ which fit into diagram (\ref{comminv}): they can be obtained one from each other by composition with the covering involution $ i_A $, i.e. the involution which exchanges the sheets of $ f $. Since the covering involution $ i_A $ is antisymplectic, i.e. $ i_A^*\omega _{X_A}=-\omega _{X_A} $, we deduce that one involution will be symplectic and the other antisymplectic. In all what follows we will denote by $ \hat{\imath} $ the unique symplectic involution induced on $ X_A $ by $ i $.
\end{rem}

In order to apply Theorem \ref{syminv} to $ \hat{\imath} $ we need to show that there are Lagrangian subspaces $A\in \mathbb{LG}^0(\wedge ^3 V)$ invariant for $i$ such that $ X_A $ is smooth.

Given the decomposition $V=V^+\oplus V^-$ as direct sum of eigenspaces of $  i$, we get 
\begin{equation}\label{v3es}
\wedge ^3V=(\wedge ^3V)^+\oplus(\wedge ^3V)^-=\wedge ^3V^+\oplus (V^+\otimes\wedge ^2V^-)\oplus (V^-\otimes\wedge ^2V^+)\oplus \wedge ^3V^-
\end{equation}
A subspace $ A\subset \wedge ^3V$ is invariant under $ i $ if and only if it can be written $ A^+\oplus A^- $, with $ A^+\subset \wedge ^3V^+\oplus (V^+\otimes\wedge ^2V^-) $ and $ A^-\subset (V^-\otimes\wedge ^2V^+)\oplus \wedge ^3V^- $.

We need to check that for such a general Lagrangian $ A $ we have $ X_A $ smooth, i.e. that $ A $ does not contain any decomposable tensors and that $ \dim A\cap F_l\leq 2 $ for all $ l\in \pn (V) $ (see \cite{OG}).

\begin{pro}
If $  \dim V^+=5 $ or $ 3 $, $ X_A $ is not smooth.
\end{pro}
\textbf{Proof.}
If $ \dim V^+=5 $, $ X_A $ is not smooth.
Indeed, we have either $ \dim A^+\geq 5 $ or $ \dim A^-\geq 5 $. In the first case, since in $ \pn (\wedge ^3V^+)\cong\pn^9 $ decomposable tensors are parametrized by $ G(2,5) $ of dimension 6, it must be $ \pn(A^+)\cap G(2,5)\neq\emptyset $. If $ \dim A^-\geq 5 $, since $ \pn (V^-\otimes\wedge ^2 V^+)\cong \pn^9 $ and decomposable tensors are parametrized by $\pn(V^-)\times G(2,5) $ of dimension 6, we get $ \pn(A^-)\cap (\pn(V^-)\times G(2,5))\neq\emptyset $.

If $\dim V^+=3$, $ X_A $ is not smooth. Indeed, if $ \dim A^+\geq 6 $ we immediately find a decomposable tensor in $ \pn(A^+) $, since $ \pn (V^+)\times G(2,3) $ is a subvariety of dimension 4 of the variety of decomposable tensors in $ \pn^9 $. An analogous dimensional count shows that there is a decomposable in $ \pn (A^-) $ if $ \dim A^-\geq 6 $.

Let us suppose that $ \dim A^+=\dim A^-=5 $. We claim that there is $ v\in V^+ $ such that $\dim (v\wedge\wedge ^2 V)\cap A\geq 3$ and this shows that $ A $ is not in $ \mathbb{LG}(\wedge ^3V)^0 $.

First of all, let us remark that there is $ v\in V^+ $ such that $(v\wedge\wedge ^2 V)\cap A^+\neq 0$  if and only if $ p(A^+)\cap (v\wedge\wedge^2V^-)\neq 0 $ where $ p:\wedge^3V^+\oplus (V^+\otimes\wedge^2V^-)\longrightarrow V^+\otimes\wedge^2V^- $ is the projection. Indeed, either $ A^+\cap \wedge^3 V^+\neq 0 $, and in this case there is a decomposable in $ A^+ $ and the proof ends, or $ \dim p(A^+)=5$. Hence there is a finite number of vectors $ v\in V^+ $ such that $(v\wedge\wedge ^2 V)\cap A^+\neq 0$: indeed, we have seen that $ \dim \pn (p(A^+))=4 $ and on the other hand $\pn(V^+)\times\pn(\wedge ^2V^-)  $ is a 4-dimensional subvariety of $ \pn( V^+\otimes\wedge^2V^-)\cong \pn^8$, so they intersect in a finite number of points.

Let $ v\in V^+ $ such a vector; let us show that $ \dim (v\wedge\wedge ^2 V)\cap A^-\geq 2$. Since $ A^-\subset (\wedge^2V^+\otimes V^-)\oplus\wedge^3V^-$, this is equivalent to $ \dim (v\wedge V^+\otimes V^-)\cap A^-\geq 2$. Let us write $ A':=A^-\cap (\wedge^2V^+\otimes V^-) $; we have $ \dim A'=4 $ otherwise there would be a decomposable in $ A^- $. Hence the morphism $ \wedge v:A'\longrightarrow\wedge^3V^+\otimes V^- $ has nontrivial kernel $ W $. 

On the other hand we have shown that there is $ \tau\in\wedge^2V\setminus \lbrace 0\rbrace $ such that $ v\wedge \tau\in A^+ $ is nonzero. The composition $ \wedge\tau\circ \wedge v:A'\longrightarrow\wedge^3V^+\otimes V^-\longrightarrow\wedge^6V $ must then be zero, since by hypothesis $ A $ is Lagrangian and this happens only if $ A^+\wedge A^-=0 $. It follows that $ \dim W\geq 2 $, since otherwise $ \wedge v $ would be surjective and $ \wedge \tau $ should be identically zero, which would imply $ \tau\in\wedge^2V^+ $ and give a decomposable tensor in $ A $.
\qed\\


We are left with the case in which $ \dim V^+=4 $, but we need a deeper analysis to understand it. First of all let us remark that in this case $ i $ is symplectic on $ \wedge^3V $, since $ \det i=1 $, and this implies also that the decomposition (\ref{v3es}) of $ \wedge^3V $ in eigenspaces of $ i $ is orthogonal with respect to the symplectic form $ \omega $: if $ \alpha\in(\wedge ^3V)^+ $ and $ \beta\in(\wedge ^3V)^- $, $\det i(\alpha \wedge \beta)=i( \alpha) \wedge i(\beta )=-\alpha \wedge \beta$, hence $ \alpha \wedge \beta=0 $.

Let us also recall a standard fact from linear algebra.
\begin{rem}\label{linalg}
If we have 3 vector spaces $ W $, $ E_1$ and $E_2 $ such that:
\begin{enumerate}
\item $ \dim W=\dim E_1=\dim E_2 $;
\item $ W\subset E_1\oplus E_2 $;
\item $ W\cap E_i=0 $ for $ i=1,2 $,
\end{enumerate}
then there is an isomorphism $ f:E_1\longrightarrow E_2 $ such that $ W $ is the graph of $ f $.

Since the third assumption implies that the two projections $ p_i: W \longrightarrow E_i$ are isomorphisms, the isomorphism $ f=p_2\circ p_1^{-1} $ makes the deal.
\end{rem}

\begin{lem}\label{a+-}
If $ \dim V^+=4 $ and $ f_1,f_2 $ are a basis of  $ V^- $, let $ A=A^+\oplus A^- $ be a Lagrangian subspace of $\wedge ^3V$ such that $ \wedge^3V^+\cap A^+=0 $, $ A^+\cap (\wedge^2V^-\otimes V^+)=0 $ and $ A^-\cap (\mathbb{C}f_i\otimes \wedge ^2 V^+)=0 $ for $ i=1,2 $. Then: 
\begin{enumerate}
\item there is a self-adjoint operator $u:\wedge ^2V^+\rightarrow \wedge ^2V^+$ such that $A^-=\left\lbrace f_1\wedge x+f_2\wedge u(x) /x\in\wedge ^2V^+\right\rbrace $;
\item $ A^+=\left\lbrace f_1\wedge f_2\wedge v+\phi(v)/v\in V^+ \right\rbrace  $ where $ \phi:V^+\rightarrow \wedge ^3V^+ $ is a linear isomorphism such that $ v\wedge\phi(w)=w\wedge\phi(v) $ for all $ v,w\in V^+ $.
\end{enumerate}
\end{lem}
\textbf{Proof.}  $ A $ is Lagrangian if for all $ v,w\in A $ we have $ v\wedge w=0 $. Here $ A $ will be Lagrangian as soon as $A^+$ and $A^-$ are Lagrangian respectively in $\wedge ^3V^+\oplus (V^+\otimes\wedge ^2V^-)$ and in $(V^-\otimes\wedge ^2V^+)$, since $A^+\wedge A^-=0 $ comes from the orthogonality of the decomposition (\ref{v3es}).

\begin{enumerate}
\item First of all let us remark that given $ u $ and $ A^- $ as in the statement, $ A^- $ is Lagrangian: for all $ x,y\in\wedge^2V^+ $ we have 
\begin{equation}\label{symm}
(  f_1\wedge x+f_2\wedge u(x))\wedge( f_1\wedge y+f_2\wedge u(y))=f_1\wedge f_2\wedge (-x\wedge u(y)+u(x)\wedge y)=0\end{equation}
Now let us consider $ A^-\subset V^-\otimes \wedge ^2 V^+=(\mathbb{C}f_1\otimes \wedge ^2 V^+)\oplus (\mathbb{C}f_2\otimes \wedge ^2 V^+) $;
then by Remark \ref{linalg} there is an isomorphism
$u:\wedge ^2V^+\rightarrow \wedge ^2V^+$ such that $ A^- $ is its graph; (\ref{symm}) tells us that $ u $ is self-adjoint because $ A^- $ is a Lagrangian subspace.

\item If $ A^+ $ and $ \phi $ satisfies the statement, $ A^+ $ is Lagrangian in $ \wedge^3V^+\oplus(\wedge^2V^-\otimes V^+ )$: indeed, for all $ v,w \in V^+ $ we have
\begin{equation}\label{bili}
  (f_1\wedge f_2\wedge v+\phi(v))\wedge(f_1\wedge f_2\wedge w+\phi(w))= f_1\wedge f_2\wedge (v\wedge\phi(w)+ \phi(v)\wedge w)=0\end{equation}
Viceversa, given a Lagrangian subspace $ A^+ $ in $\wedge^3V^+\oplus(\wedge^2V^-\otimes V^+ )$ such that $ \wedge^3V^+\cap A^+=0 $ and $ A^+\cap (\wedge^2V^-\otimes V^+)=0 $, 
by Remark \ref{linalg} there is an isomorphism
$ \phi:V^+\rightarrow \wedge ^3V^+ $ linear such that $ A^+ $ is its graph. Since $ A^+ $ is Lagrangian, from (\ref{bili}) we deduce that $ \phi $ satisfies $ v\wedge\phi(w)=w\wedge\phi(v) $ for all $ v,w\in V^+ $.\qed
\end{enumerate}

\begin{lem}\label{a-}
Using the notation of Lemma \ref{a+-}, if $ u $ has 6 distinct eigenvalues and no decomposable eigenvector in $ \wedge^2V^+ $, then $ A ^-$ does not contain any decomposable tensor and there is a basis of eigenvectors $ x_1,\dots,x_6 \in\wedge^2V^+$ such that $ u $ is diagonalizable.
\end{lem}
\textbf{Proof.} If $ v\wedge w_1\wedge w_2\in A^- $ we can suppose $v\in V^- $ and it follows that there must be $ \lambda\in\mathbb{C} $ such that $u( w_1\wedge w_2)=\lambda w_1\wedge w_2 $. This is against our assumption that $ u $ has no decomposable eigenvectors, hence there are no decomposable tensors in $ A^- $. Let $ x\in \wedge^2V^+ $ be an eigenvector of $ u $; since $ Q(x)\neq 0 $, we have an orthogonal decomposition $ \wedge^2V^+=\mathbb{C}x\oplus V' $ such that $ V' $ is invariant for $ u $. By iterating this reasoning we then get an orthogonal basis of eigenvectors.\qed

\begin{lem}\label{a+}
Using the notation of Lemma \ref{a+-}, $ A ^+$ does not contain any decomposable tensor.
\end{lem}
\textbf{Proof.} If $ v\wedge w_1\wedge w_2\in A^+ $ then we can suppose $v\in V^+ $ and it follows that $ w_1\wedge w_2 $ is a decomposable in $ \wedge^2 V^+\oplus\wedge^2 V^- $ which is only possible if $ w_1\wedge w_2 \in \wedge^2 V^+ $ or $ w_1\wedge w_2 \in \wedge^2 V^- $, against the fact that  $ A^+\cap \wedge^3V^+=0 $ and $ A^+\cap (\wedge^2V^-\otimes V^+)=0 $.\qed\\

Let us define $ \mathbb{LG}(\wedge^3V)^*$ to be the set of all $ A\in \mathbb{LG}(\wedge^3V)$ such that $A$ admits a decomposition as in Lemma \ref{a+-} with $ u $ satisfying also the hypothesis of Lemma \ref{a-}. It follows from what we said that it is an open set inside $  \mathbb{LG}(\wedge^3V) $. Indeed, up to a base change of $ V^- $, $ \wedge^3V^+\cap A^+=0 $, $ A^+\cap (\wedge^2V^-\otimes V^+)=0 $ and $ A^-\cap (\mathbb{C}f_i\otimes \wedge ^2 V^+)=0 $ for $ i=1,2 $ are all open conditions, since $ \pn(A^-) $ and $ \pn(\mathbb{C}f_i\otimes \wedge ^2 V^+) $ are both 5-dimensional linear subspaces of $ \pn^{11} $ and $\dim\pn(A^+)= \dim\pn( \wedge^3V^+)=\dim\pn( \wedge^2V^-\otimes V^+)=3$ in $ \pn^7 $.


\begin{pro}\label{nodec}
If $ A=A^+\oplus A^-\in \mathbb{LG}(\wedge^3V)^* $ then $ A\in \mathbb{LG}(\wedge ^3V)^0 $.
\end{pro}
\textbf{Proof.} Let us remark that if $ v\wedge w_1\wedge w_2\in A $, such a decomposable is associated to a 3-dimensional vector subspace $ W\subset V $ which must therefore verify $\dim W\cap V^+\geq 1 $. We can hence suppose that $ v\in \pn(V^+)\cap Y_A $ and hence $ v $ is fixed by the involution $ i $. Moreover it must be $ \dim (v\wedge \wedge^2V)\cap A\geq 2 $, since otherwise we would have found a decomposable inside $ A^+ $ or $ A^- $, which is not possible by Lemma \ref{a-} and \ref{a+} since $A\in \mathbb{LG}(\wedge^3V)^* $.\\

To conclude the proof we need to analyze better the fixed points of $ i $ on $ Y_A $ and for this purpose we need to recall here the construction of the quadric line complex: for further details and all the proofs the reader is referred to Chapter 6 of \cite{GH}.

Lines $ l_x $ in $ \pn^3 $ are parametrized by the Grassmannian $ G=G(2,4)\subset\pn^5 $, which is defined by the Pl\"{u}cker quadric equation $ x\wedge x =0 $. Given another smooth quadric $ F $ in $ \pn^5 $, the intersection $ X=F\cap G $ is the so-called quadric line complex. Given $ p\in\pn^3 $ we want to understand which lines of our complex pass through $ p $. Let $ \sigma(p) $ be the set of all $ x\in X $ such that $ p\in l_x $; it is a 2-plane contained in $ G $ and we look at its intersection with $ F $, which is a conic in $ \sigma(p) $. The set $ S $ of points $ p\in \pn^3 $ such that $ F\cap \sigma(p) $ is a singular conic is shown to be a singular Kummer surface of degree 4, called the associated Kummer surface of $ X $. The singular locus $ R $ of $ S $ is made of 16 ordinary double points, which are precisely the ones such that $ F\cap \sigma(p) $ is a double line.

Given $ x\in X $, the corresponding line $ l_x $ is singular if there is $ p \in l_x $ such that $ \sigma(p) $ is tangent to $ F $ in $ x $. The set $ \Sigma $ of points $ x\in X $ such that $ l_x $ is singular is a smooth minimal K3 surface and there is a morphism $ \pi:\Sigma \longrightarrow S $ defined by $ \pi(x)=p\in l_x $ such that $ \sigma(p)=T_xF $. In fact, $ \pi $ is the blow-up of $ S $ along $ R $. There is also a morphism $ \pi':\Sigma \longrightarrow S^* $ defined by $ \pi'(x)=h\supset l_x $ such that $ \sigma(h)=\{y\in G/l_y\subset h\} $ is tangent to $ F $ in $ x $. There is a commutative diagram
\[ 
\xymatrix{
\Sigma\ar@{->}[rr]^{\pi}\ar@{->}[dr]^{\pi'}&&S\ar@{->}[dl]^{\delta}\\
&S^*&} \]
where $ \delta(v)=T_vS $ for all $ v\in S $. 

Another characterization of $ \Sigma $ is the following. Let  $x\wedge Qx =0$ and $ x\wedge Q'x=0 $ be the two quadratic forms defining the two quadrics $ G $ and $ F $, where $Q$ and $ Q'$ are two symmetric matrices; $ \Sigma=G\cap F\cap H $ where $ H $ is the quadric hypersurface corresponding to the matrix $ Q'Q^{-1}Q' $. Since we will need it later, let us remark that, by standard linear algebra, whenever $Q'$ has distinct eigenvalues it is possible to suppose that $ Q =I$ and at the same time to diagonalize $ Q' $ (see Chapter XII \S 6 of \cite{Ga}): hence we can find homogeneous coordinates $ [X_0,\dots,X_5]\in \pn^5$ such that $ G $ and $ F $ are respectively defined by the equations $ \sum_{i=0}^5 X_i^2=0 $ and $ \sum_{i=0}^5 \lambda_i X_i^2=0 $. Using such coordinates then $ H $ turns out to be defined by $ \sum_{i=0}^5 \lambda_i^2 X_i^2=0 $.

\begin{pro}\label{prfix}
Let $ i $ be an involution on $ V $ such that $ \dim V^+=4 $ and let $ A\in \mathbb{LG}(\wedge ^3V)^* $ be invariant for $ i $. Then the fixed locus $ \Fix(i) $ of $i$ on $Y_A$ is the union of 6 isolated fixed points $ q_1,\dots ,q_6 \in \pn(V^+)$, one smooth quadric $ Q $ and a singular Kummer surface $ S $ of degree $ 4 $ in $ \pn ^3 $.
\end{pro}

\textbf{Proof.} The fixed points of $ i $ on $ Y_A $ are precisely the intersections $Y_A\cap \pn(V^+) $ and $  Y_A\cap \pn(V^-) $. Given $ v\in V $ we want to understand when $ (v\wedge \wedge ^2V) \cap A\neq 0 $; since $ (v\wedge \wedge ^2V) \cap A $ is fixed by the involution, this intersection splits into \[(v\wedge \wedge ^2V) \cap A =((v\wedge \wedge ^2V)^+ \cap A ^+)\oplus((v\wedge \wedge ^2V)^-\cap A ^-) \]where $ (v\wedge \wedge ^2V)^+  $ and $ (v\wedge \wedge ^2V)^- $ are the intersections of $(v\wedge \wedge ^2V) $  respectively with $ (\wedge^3V)^+ $ and $ (\wedge^3V)^- $. We are going to investigate each of these summands separately.

\begin{enumerate}
\item If $ v\in V^+ $ we have $ (v\wedge \wedge ^2V)^+= v\wedge( \wedge ^2V^+\oplus \wedge ^2V^- ) $ and $  (v\wedge \wedge ^2V)^-=v\wedge (V^+\otimes V^-) $.
\begin{itemize}
\item If $ \alpha \in (v\wedge \wedge ^2V)^+ \cap A ^+ $, there is $ \tau\in\wedge^2V^+ $ such that $ \alpha=v\wedge \tau+ v\wedge f_1\wedge f_2$ and on the other hand since $ \alpha\in A^+ $ we have $ \alpha= f_1\wedge f_2\wedge v+\phi(v) $ by Lemma \ref{a+-}. Hence we get $ \phi(v)=v\wedge \tau $ and this happens if and only if $ v\wedge\phi(v)=0 $: this equation defines a quadric $ Q $ in $ \pn^3\cong\pn (V^+) $.
\item If $ \alpha\in (v\wedge \wedge ^2V)^-\cap A^- $, there are $ y_1,y_2\in V^+ $ such that $ \alpha=v\wedge y_1\wedge f_1+v\wedge y_2\wedge f_2 $ and on the other hand there is $ x\in\wedge^2V^+ $ such that $ \alpha=f_1\wedge x+f_2\wedge u(x)$ by Lemma \ref{a+-}. Comparing the two expressions we get
\begin{equation}\label{sysxux}
\left\lbrace 
\begin{array}{l}
x=v\wedge y_1 \\
u(x)=v\wedge y_2
\end{array}\right.  
\end{equation} 
and the forms $ x\in\wedge^2V^+ $ solutions of (\ref{sysxux}) are exactly those who satisfy $ x\wedge x=x\wedge u(x)=u(x)\wedge u(x)=0 $ in $ \pn(\wedge^2 V^+)\cong \pn^5 $. Hence $ x\in\Sigma $, the smooth K3 surface associated to the quadric line complex defined above, where as $ F $ we consider the quadric hypersurface defined by $ x\wedge u(x)=0 $.

Keeping the notation we used above, we claim that $ \pi(x)=v\in S $. Indeed, $ \pi(x)\in l_x $ is the point such that $ T_xF=\sigma(\pi(x)) $; hence every line corresponding to a point of $ T_xF $ passes through $ \pi(x) $, which can be recovered as the intersection of any such line with $ l_x $. Since $ T_xF $ is defined by the equation $ y\wedge u(x)=0 $ for $ y\in\pn^5 $, we get that $ u(x)\in T_xF $, so $ \pi(x) $ is the intersection of the lines $ l_x $ and $ l_{u(x)} $, i.e. $ \pi(x)=v $.



\end{itemize}
\item If $ v\in V^-$ we have $ (v\wedge \wedge ^2V)^+=  v\wedge (V^+\otimes V^-) $ and $  (v\wedge \wedge ^2V)^-=v\wedge \wedge ^2V^+ $.
\begin{itemize}
\item We have $A^+\cap (V^+\otimes \wedge ^2V^-)=0$ by Lemma \ref{a+-}, so there are no fixed points which arise from this case.
\item If $ \alpha\in (v\wedge \wedge ^2V)^-\cap A^- $, there is $ \tau\in \wedge ^2V^+ $ such that $  \alpha=v\wedge \tau $ and on the other hand there is $ x\in\wedge^2V^+ $ such that $ \alpha=f_1\wedge x+f_2\wedge u(x)$ by Lemma \ref{a+-}. Since $f_1,f_2$ form a basis of $V^-$, there is $ \lambda\in\mathbb{C} $ such that $ v=f_1+\lambda f_2 $. Comparing the two expressions we get $ \tau=x $ and $ u(x)=\lambda x $, i.e. 
 $\lambda $ is an eigenvalue of $ u $ and $ x $ is the corresponding eigenvector. Since $ u $ has 6 different eigenvalues, we obtain 6 isolated fixed points.\qed

\end{itemize}
\end{enumerate}

\begin{rem}\label{dim2}
When $ v\in S $ the proof of Proposition \ref{prfix} gives us more information: indeed,  $ \pi^{-1}(v)\cong\pn((v\wedge \wedge ^2V)^-\cap A ^-) $, hence $ \dim(v\wedge \wedge ^2V)^-\cap A^-=2 $ if $v\in R  $, $ 1 $ otherwise.

When $ v\in Q $ we always have that $ \dim(v\wedge \wedge ^2V)^+\cap A^+=1 $, since by Lemma \ref{a+-} the projection $ A^+\rightarrow V^+\otimes \wedge^2V^- $ is an isomorphism.
\end{rem}

\textbf{End of the proof of Proposition \ref{nodec}.} From Remark \ref{dim2} it follows that the only points $v$ for which $ \dim (v\wedge \wedge^2V)\cap A= 3 $ are among the 16 isolated singular points of $ S $ if they belong to $ Q $ too,
which does not happen if $ Q $ is general enough. Hence we have $ \dim (v\wedge \wedge^2V)\cap A\leq 2 $.

Moreover, we have seen that if $ v\wedge w_1\wedge w_2\in A $ we can suppose $ v\in V^+ $ and we must have $ \dim (v\wedge \wedge^2V)\cap A= 2 $. In such a case then $ v\in Q\cap S $; we claim that $ v $ is then a singular point of the intersection, which cannot happen when $ Q$ is general enough.

Indeed, there is $ x\in\wedge^2V^+ $ such that $  v\wedge w_1\wedge w_2= v\wedge f_1\wedge f_2+\phi(v)+f_1\wedge x+ f_2\wedge u(x) $ and by the proof of Proposition \ref{prfix} we know that there are $ z\in \wedge ^2V^+ $ and $ y_i\in V^+ $ such that $ \phi(v)=v\wedge z $, $ x=v\wedge y_1 $ and $ u(x)=v\wedge y_2 $. We can suppose that $ w_i=w_i^++f_i $ with $w_i^+\in V^+  $ for $ i=1,2 $ up to a base change: we cannot have $ W\subset V^+\oplus \mathbb{C}f_i $ for $ i=1 $ or $ 2 $, because in that case $ v\wedge w_1 \wedge w_2\in \wedge^3 V^+\oplus(\mathbb{C}f_i\otimes\wedge^2V^+) $ and $ (\wedge^3V^+\oplus(\mathbb{C}f_i\otimes\wedge^2 V^+))\cap A=0 $. After replacing we get
\begin{equation}\label{sysform}
\left\lbrace \begin{array}{l}
v\wedge w_1^+\wedge w_2^+=v\wedge z\\
 \\
v\wedge w_1^+\wedge f_2-v\wedge w_2^+\wedge f_1=v \wedge y_1\wedge f_1+v \wedge y_2\wedge f_2
\end{array}
\right.
\end{equation}

From the second equation we get $ w_1^+- y_2=kv $ and $ w_2^++ y_1=hv $ with $ k,h\in \mathbb{C} $. Hence we have $ w_1^+\wedge w_2^+=y_1\wedge y_2+hy_1\wedge v-ky_2\wedge v  $. Replacing in the first equation we then obtain $ v\wedge (y_1\wedge y_2-z)=0 $. This shows that $ y_1,y_2\in T_vQ $, which is defined by $ y\wedge \varphi (v)=0 $.


On the other hand, $ T_vS $ is spanned by the two lines associated to $ x=v\wedge y_1 $ and $ u(x)=v\wedge y_2 $: 
both $ x $ and $ u(x) $ satisfy $ y\wedge u(x) =0$ which is the equation defining $ T_xF $, hence they span $ \pi'(x)=T_vS $.
Thus we get $T_vS=T_vQ $,
 i.e. $ v $ is a singular point of $ Q\cap S $. This ends the proof of the smoothness of $ X_A $ when $A\in \mathbb{LG}(\wedge^3V)^* $.\qed\\

Let us now study the symplectic involution $ \hat{\imath} $ induced on $ X_A $ (see Remark \ref{indinv}).


\begin{pro}
Let $ i $ be an involution on $ V $ such that $ \dim V^+=4 $ and let $ A\in \mathbb{LG}(\wedge ^3V)^* $ be invariant for $ i $; then the induced symplectic involution $ \hat{\imath} $ on $ X_A $ has 28 isolated fixed points and a fixed K3 surface.
\end{pro}
\textbf{Proof.} We know that $ \Fix(\hat{\imath}) $ has smooth symplectic components from Lemma \ref{lsmo} and Proposition \ref{psymsub} and on the other hand we have $ f (\Fix(\hat{\imath}))\subset \Fix(i) $ which we have completely described.
 
If $ Z\subset \Fix(\hat{\imath}) $ is a surface it must be the double cover either of $ Q $ or of $ S $ and we know which is the ramification locus: it is given in the former case by $ Q\cap W_A $ and in the latter by $ S\cap W_A $.

From what we said above it is clear that $ Q\cap W_A\cong Q\cap S $ and the double cover of a smooth quadric ramified along a quartic curve is a K3 surface. On the other hand let $ C $ be the trace of the quadric $ Q $ on the Kummer surface $ S $; then $ S\cap W_A $ is the union of the 16 ordinary double points of $ S $ and of $ C=Q\cap S $. In this case we have the following commutative diagram
$$
\xymatrix{
T\ar@{->}[r]^{\varepsilon}\ar@{->}[d]^g&Z \ar@{->}[d]^{f_{|Z}}\\
\Sigma \ar@{->}[r]^{\pi}& S}$$
where $ \pi $ and $ \varepsilon $ are respectively the blow-ups of $ S $ and $ Z $ in $ p_1,\ldots ,p_{16} $ and in $ f^{-1}(p_1),\dots,f^{-1}(p_{16}) $ and $ g $ is the double cover of $ \Sigma $ ramified along $ E_1,\dots,E_{16} $ and $ \pi ^{-1}(C)\cong C $. Let $ D $ be a divisor on $ \Sigma $ such that $ 2D= C$ and $ F_i $ be the exceptional divisor on $ T $ corresponding to $ f^{-1}(p_i) $; we have that $ K_T=g^*D+\sum F_i $ and also $ K_T=\varepsilon^*K_Z+\sum F_i $,  hence $ K_Z$ cannot be trivial and $ Z $ is neither K3 or abelian.


From what we said it follows that there cannot be abelian surfaces inside $\Fix(\hat{\imath})$ and that the symplectic involution is of the third type described by Theorem \ref{syminv} and it must fix $ 28 $ isolated points and a K3 surface. We have already seen that the K3 surface arises as the double cover of $ Q $ ramified along the quartic curve $ Q\cap S $. Moreover the 6 isolated fixed points on $ Y_A $ gives 12 fixed points for $ \hat{\imath} $. Finally the 16 points $ f^{-1}(p_i) $ which are the fibers of the 16 ordinary double points of $ S $ are fixed too, giving us all the 28 isolated points we expected.\qed


\bibliographystyle{plain}
\bibliography{th_bib}

\end{document}